\documentclass[12pt]{article}
\usepackage{amssymb}

\newtheorem{thm}     {Theorem}[section]

\newtheorem{definition}  [thm]{Definition}
\newtheorem{cor}     [thm]{Corollary}
\newtheorem{lemma}   [thm]{Lemma}

\newcommand{\proof} {\noindent{\bf Proof. }}

\newcommand{\B}{\mathbb B}
\newcommand{\C}{\mathbb C}
\newcommand{\D}{\mathbb D}

\newcommand{\R}{\mathbb R}

\newcommand{\HH}{\mathbb H}

\newcommand{\OO}{\mathcal O}

\def\Re{{\rm Re\,}}
\def\Im{{\rm Im\,}}

\begin{document}

\title{On   the Wong-Rosay  theorem}
\author{Alexandre Sukhov }
\date{}
\maketitle

\bigskip
 
{\small Abstract.  We prove a local version of the Wong - Rosay theorem for a domain with a piecewise smooth generic 
strictly pseudoconvex boundary point.}



\bigskip

\section{Introduction}

Let $\Omega$ be a domain with non-empty boundary $b\Omega$ in a complex manifold $M$ of complex dimension $n$.   
\begin{definition}
\label{point}
 We say that $p \in b\Omega$ is a {\it piecewise smooth generic strictly pseudoconvex (boundary) point} if the following assumptions hold:
\begin{itemize}
\item[(i)] there exists an open  connected neighborhood $U$ of $p$ in $M$ such that 
\begin{eqnarray}
\label{BoundPoint}
\Omega \cap U = \{ z \in U: \rho_j(z) < 0, j=1,...,m\}
\end{eqnarray}
Here every function $\rho_j$ is $C^2$  strictly plurisubharmonic on $U$, and $\rho_j(p) = 0$.  
\item[(ii)] $\partial \rho_1 \wedge ... \wedge \partial \rho_m \neq 0$ on $U$.
\end{itemize} 
\end{definition}
  Of course, the condition (i) can be stated in the equivalent form: the hypersufaces $\Gamma_j = \{ \rho_j = 0 \}$ 
  ({ \it the local faces} of $b\Omega$) are strictly pseudoconvex i.e. the Levi from of each $\Gamma_j$ is   positive defined on the holomorphic tangent bundle of $\Gamma_j$.
  The condition (ii) assures that the real submanifold $\{ \rho_j = 0, j=1,...,m \}$  ({\it the corner}) is generic.
  A point $p$ is {\it smooth} if $m=1$. Of course, in this case  a boundary point $p$ is a usual $C^2$ strictly pseudoconvex point.
 
  We  denote by $Aut(\Omega)$
 the holomorphic automorphism group of $\Omega$ equipped with the compact open topology. 
{\it The limit set of} $Aut(\Omega)$ is the set of points $p \in b\Omega$ such that there exists a point $q \in \Omega$
and a sequence of automorphisms $(f^k)_k$ in $Aut(\Omega)$ satisfying  $\lim_{k \to \infty} f^k(q) = p$.

If $z = (z_1,...,z_n)$ are the standard coordinates of $\C^n$, we write $z = (z_1,z')$ with $z' = (z_2,...,z_n)$, and $z_j = x_j + iy_j$ with $x_j,y_, \in \R$. Also, $\parallel z \parallel^2 = \sum \vert z_j \vert^2$ denotes the Euclidean norm. In what follows we use the notation $$\B^n = \{ z \in \C: \parallel z \parallel^2 < 1\}$$ for the Euclidean unit ball of $\C^n$  and we denote by 
\begin{eqnarray*}
\label{Siegel}
\HH = \{ z \in \C^n: \Re z_1 + \parallel z' \parallel^2 < 0\}
\end{eqnarray*}
the unbounded realization of $\B^n$ (recall that the domain $\HH$ is biholomorphic to $\B^n$ via the Caley transform).

The goal of the present paper is to prove   the following

\begin{thm}
\label{theo1}
Assume that the limit set of $Aut(\Omega)$ contains  a piecewise smooth generic strictly 
pseudoconvex  point $p \in b\Omega$.  
Then  $\Omega$ is biholomorphic to $\B^n$ and $p$ is a smooth strictly pseudoconvex point.
\end{thm}
This result is definitive: neither the condition (i) nor the condition (ii) in Definition \ref{point} can be dropped as show the following examples.

First,  consider  the domain 
$$\Omega_1 = \{ \rho_1 =  \Re z_1 + \vert z_2 \vert^2 < 0,  \rho_2 = \Re z_2 < 0 \}$$ 
which is  invariant with respect to the 1-parameter family of dilations 
$d_t: (z_1, z_2) \mapsto (t z_1, \sqrt{t} z_2)$, $t > 0$. This family is non-compact because it degenerates
when $t = 0$ and the origin belongs to the limit set of $Aut(\Omega_1)$.  However, $\Omega_1$ is not biholomorphic to $\B^2$. The domain  $\Omega_1$ satisfies (ii), but does not satisfy (i): one of the faces is not strictly pseudoconvex.

Second, consider the domain
$$\Omega_2 = \{ z \in \C^2: \rho_1 = \Re z_1 + \vert z_2 \vert^2 < 0, \rho_2 = \Im z_1 + \vert z_2 \vert^2 < 0 \}$$
invariant with respect to the same family of dilations $(d_t)$. Of course, $\Omega_2$ is not biholomorphic to $\B^2$ as well. Here the assumption (i) holds, but (ii) fails at the origin (thought $d \rho_1 \wedge d \rho_2 \neq 0$).

Theorem \ref{theo1} belongs to a series of results which often are called { \it the Wong - Rosay type theorems}. This short paper is not an appropriate place in order to present 
a detailed history and the state of art of this direction. Hence I restrict myself by results directly concerning the topic of present  paper; in particular, I skip a discussion of (many beautiful) results dealing with the non-strictly pseudoconvex  case.

The fact that a bounded strictly pseudoconvex domain in $\C^2$ with the maximal possible dimension (equal to $8$) of $Aut(\Omega)$ (which a real Lie group ) is biholomorphic to 
$\B^2$, was  known already to Elie Cartan \cite{Ca}. One can view this phenomenon as a special case of the general differential geometric principle stating that structures with 
rich automorphism groups usually are flat.  In \cite{BS} Burns - Shnider proved that a bounded strictly pseudoconvex domain $\Omega$ in $ \C^n$ with non-compact automorphism group is biholomorphic to the unit ball. This was a striking and surprising result because the assumption of non-compactness of the group $Aut(\Omega)$ {\it a priori } does not impose restrictions on the  dimension of $Aut(\Omega)$. Their proof is based on the Chern - Moser theory \cite{CM} (more precisely they use the part due to Chern which extends the approach of E. Cartan to higher dimension; the approach of Moser is very different ) and requires a relatively high regularity (at least, $C^6$) of $b\Omega$. The group $Aut(\Omega)$ is non-compact if and only if its limit set on the boundary of $\Omega$ is not empty. Wong \cite{W} and Rosay \cite{Ros} discovered that the result of Burns - Shnider can be localized: it suffices to assume that the limit set of $Aut(\Omega)$ contains a strictly pseudoconvex point under the assumption that $\Omega$ is bounded. Perhaps, their most important observation was that the phenomenon discovered by Burns - Shnider , in fact, can be treated without the Cartan - Chern - Moser approach. It turned out that other geometric tools (such as biholomorphically invariant metrics and normal families of holomorphic maps)  are  more efficient and lead to more general results.  Later their approach was considerably simplified by Pinchuk \cite{Pi} using his version of so called scaling method. The first  purely local version of this phenomenon  was  obtained by Efimov \cite{Ef}; he established Theorem \ref{theo1} in the special case where $p$ is a smoooth strictly pseudoconvex point (i.e. $m = 1$ in Definition \ref{point}). His result was  extended  by Gaussier - Kim - Krantz \cite{GKK} to the case where $\Omega$ is a domain in a complex manifold. In the non-smooth case, Coupet - Sukhov \cite{CoSu} proved that a bounded piecewise smooth strictly pseudoconvex domain 
with generic corners in $\C^n$ is biholomorphic to the unit ball if $Aut(\Omega)$ is non-compact (and, therefore, the boundary of $\Omega$ necessarily is smooth). 

Theorem \ref{theo1} generalizes all above mentioned results beginning by the works of Wong and Rosay. Our proof consists of two parts. The first one concerns a localization of the Kobayashi - Royden metric. The second (and the principal ) part is based on the scaling method. Notice that the proof of Efimov is based on the version of this method due to Pinchuk \cite{Pi}. In the present paper we use the approach due to Frankel \cite{F1} which seems to be more adapted to the non-smooth case. A detailed presentation of various  versions and applications of the scaling method is contained in the expository paper of Berteloot \cite{Ber}. The approach of Frankel also was used in \cite{CoSu}. In the present paper we simplify the arguments  from \cite{CoSu} reducing them to the known estimates of the Kobayashi-Royden metric. This works because we deal with a special type of boundary points while the theory of Frankel  concerns general convex domains.

 \section{The Kobayashi-Royden metric and normal families}

 This section is preliminary. For the convenience of readers, I recall several results concerning the Kobayashi - Royden metric. 
 
 Fix any Riemannian metric on $M$ and use it in order to measure the distances 
 on $M$ and norms of tangent vectors. In the case where $M = \C^n$ we always use the standard Euclidean norm and metric. In what follows we denote by $\D = \{ \zeta \in \C: \vert \zeta \vert < 1 \}$ the unit disc in $\C$ (i.e. $\B^1$). Let also  $\OO(\D,M)$ denotes the space of holomorphic maps from $\D$ to $M$; we call such maps complex or analytic discs in $M$.

 Recall that {\it the Kobayashi - Royden pseudometric} $F_M$ of $M$ is defined  on a point $p \in M$ and a tangent vector $v \in T_pM$ by
 
 \begin{eqnarray*}
 F_M(p,v) = \inf\left\{ \lambda^{-1} : \, \mbox{there exists} \, f \in \OO(\D,M) \, \mbox{such that} \, f(0) = p, \, \frac{df}{d\zeta}(0) = \lambda v, \, \lambda > 0 \right\}
 \end{eqnarray*}
 Denote by $K_M(p,q)$ the usual Kobayashi pseudodistance of $M$ between points $p,q \in M$.  According to the fundamental result of Royden \cite{R}, $F_M$ is an upper semicontinous function on the tangent bundle of $M$ and $K_M$ is the integration form of $F_M$.

  We will use the fundamental property of the Kobayashi-Royden pseudometric and the Kobayashi pseudodistance: they are holomorphically decreasing. Namely, if $f: M \to N$ is a holomorphic map between two complex manifolds, then $F_N(f(p), df(p) v) \le F_M(p,v)$ and $K_N(f(p),f(q)) \le K_M(p,q)$. Recall also that $M$ is called  {\it hyperbolic at $p \in M$} if there exists a constant $C > 0$ such that $F_M(p,v) \ge C \parallel v \parallel$ for every tangent vector $v \in T_pM$. A manifold $M$ is called {\it locally hyperbolic}, if it is hyperbolic at every point. Also $M$ is called (Kobayashi) {\it hyperbolic} if $K_M$ is a distance that is $K_M(p,q) > 0$ when $p \neq q$; in this case it induces the usual topology of $M$.  According to \cite{R},  $M$ is hyperbolic if and only if it is locally hyperbolic. Recall also \cite{R} that $M$ is hyperbolic if and only if the family $\OO(\D,M)$ is equicontinuous (with respect to the usual topology of $M$). {\it The Kobayashi ball}  with centre at $p$ and of radius $\delta > 0$
 is defined as
 $$B_{K_M}(p,\delta) = \{ q \in M : K_M(p,q) < \delta \}$$
 Recall also that a manifold $M$ is called {\it complete hyperbolic} if it is a complete space with respect to the Kobayashi distance that is every Kobayashi ball is compactly contained in $M$.

 \subsection{Localization and normal families} Here  we discuss  some results on localization and asymptotic behavior of the Kobayashi-Royden metric. { \it Everywhere through this paper $C > 0$ denotes a positve constant which is allowed to change  its value from estimate to estimate.}

 We begin with   the following localization principle which follows from Lemma 2.2 of \cite{ChiCoSu}:

\begin{lemma}
\label{localization}
Let $p \in b\Omega$ be a picewise smooth strictly pseudoconvex point. There exist open neighborhoods $U \subset U'$ of $p$ in $M$ 
and $\delta > 0$ such that for every $q \in \Omega \cap U$ the Kobayashi ball $B_{K_\Omega}(q,\delta) $ is contained in $\Omega \cap U'$.
\end{lemma}
The hypothesis of  \cite{ChiCoSu} requires an existence of a negative plurisubharmonic function on $\Omega$ which is strictly plurisubharmonic (in the generalized sense) 
in a neighborhood of $p$. Appying to the functions $\rho_j$ from (\ref{BoundPoint}) the construction from \cite{Su}, one can extend each of these functions to a plurisubharmonic function , say $\tilde \rho_j$, globally defined and negative on $\Omega$. Then the  function $\sup_j \tilde \rho_j$ satisfies the assumptions required in \cite{ChiCoSu}.

Since the Kobayashi distance is holomorphically decreasing, we obtain the following 

\begin{cor}
\label{discs}
 There exists  $\tau = \tau(\delta) > 0$ such 
for every point $q \in \Omega \cap U$ and every holomorphic map $h:\D \to \Omega$  with $h(0) = q$, one has $h(\tau \D) \subset \Omega \cap U'$. 
\end{cor}

It follows now from the definition of the Kobayashi-Royden metric that there exists a constant $C > 0$ such that
\begin{eqnarray}
\label{loc1}
F_{\Omega \cap U'}(z,v) \le C F_{\Omega}(z,v)
\end{eqnarray}
for all $z \in \Omega \cap U$ and $v \in T_z\Omega$. We refer readers to \cite{Ber} for a detailed  discussion of related results.

As one of  the consequences of these localization  results  we have the following

\begin{lemma}
\label{convergence}
In the hypothesis of Theorem \ref{theo1}, there exists a subsequence of the sequence $(f^k)$ converging to the constant map $f \equiv p$ uniformly on 
compact subsets of $\Omega$.
\end{lemma}
Choose  coordinate neighborhood $U \subset U$ of $p$ small enough such that Corollary \ref{discs} can be applied. Let $K$ be a compact subset of $\Omega$ containing the point $q$. We claim that 
$f^k(K) \subset \Omega \cap U'$ for each $k$ big enough. Consider two finite coverings of $K$ , respectively by open coordinate neighborhoods  $V_j $ and $W_j$ , $j = 1,...,N$, such that $V_j \subset W_j \subset \Omega$, and the following holds: 
\begin{itemize}
\item[(i)] $q \in V_1$;
\item[(ii)] for every $j$ there  exists a coordinate biholomorphism  $\phi_j: W_j  \to \B^n$, such that $\phi_j(V_j) = \varepsilon_j \B^n$, where $\varepsilon_j \le  \tau$ and $\tau$ is given by Corollary \ref{discs};
\item[(iii)]  one has $q^j:= \phi_j^{-1}(0) \in V_{j-1}$, $j= 2,...,N$
\end{itemize}
For every $k$ big enough, $f^k(q) \in U$. Given unit vector $v \in \C^n$, we apply Corollary \ref{discs} to the discs $h^k: \D \to \Omega$, 
$h^k: \D \ni \zeta \mapsto f^k \circ \phi_1^{-1}(\zeta v)$. This implies that $f^k(V_1) \subset U'$. Hence, there exists a subsequence, again denoted by $(f^k)$, converging uniformly on 
$\overline V_1$ to a holomorphic map $f$. Since $f(q) = p$, by the maximum principle $f \equiv p$. Then by (iii), for $k$ big enough one has  $f^k(q^2) \in U$ and a similar argument shows that $f^k(V_2) \subset U'$. Repeating this argument for all $j$, we conclude.

\begin{cor}
\label{Hyp}
$\Omega$ is a hyperbolic  domain.
\end{cor}
Indeed, let $z^0$ be an arbitrary point of $\Omega$. Then for some $k$ big enogh $f^k(z^0) \in \Omega \cap U$. But the domain $\Omega \cap U$
is biholomorphic to a bounded domain in $\C^n$ and hence is hyperbolic. 
Therefore by (\ref{loc1}) one has
$$F_\Omega(z^0,v) = F_\Omega(f^k(z^0),df^k(z^0)v) \ge C F_{\Omega \cap U}(f^k(z^0),df^k(z^0)v) \ge C \parallel v \parallel $$
Here we used the fact that $\Omega \cap U$ is hyperbolic. Hence $\Omega$ is locally hyperbolic  and so $\Omega$ is hyperbolic.

\subsection{Estimates }

We assume that $\Omega$ satisfies assumptions of Theorem \ref{theo1}.

The following upper bound on the Kobayashi-Royden infinitesimal metric $F_\Omega$ is classical:

\begin{lemma}
\label{upperbound}
There exist a constant $C > 0$  and  a (coordinate) neighborhood $U$ of $p$ in $M$ such that for each  $z \in \Omega$ and and a tangent vector $v \in T_z\Omega$ 
one has
\begin{eqnarray*}
F_\Omega(z,v) \le C\parallel v \parallel /dist(z,b\Omega)
\end{eqnarray*}
\end{lemma}
Indeed, the ball centered at $z$ and of radius $dist(z, b\Omega)$ is contained in $\Omega$ so the estimate follows by the holomorphic decreasing property
of the Kobayashi-Royden metric.

For a bound  from below   recall some results of \cite{Su}.
\begin{lemma}
\label{belowbound}
  There exists a neighborhood $U$ of $p$ in $M$ and a constant $C > 0$ such that 
 \begin{eqnarray*}
 F_\Omega(z,v) \ge C \parallel v \parallel / dist(z,b\Omega)^{1/2}
 \end{eqnarray*}
 for every $z \in \Omega \cap U$ and $v \in T_z\Omega$.
 \end{lemma}

Finally, we need estimates of the Kobayashi-Royden metric on convex domains. Let $G \subset \C^n$ be a convex domain, $p \in G$ be a point and $v$ be a vector of $\C^n$. Consider a complex line $A$ through $p$ in the direction $v$ and denote by
$$L_G(p,v) = \sup \{ \delta > 0: \B^n(p,\delta) \cap A \subset G \}.$$ 
In other words, $L(p,v)$ is the $\sup$ of radii of discs centred at $p$ and contained in $A \cap G$. The following result is due to Graham \cite{G} and Frankel \cite{F2};
a short geometric proof is obtained by Bedford - Pinchuk \cite{BP}.

\begin{lemma}
\label{BedPin}
Let $G$ be a convex domain in $\C^n$. For every $p \in \Omega$ and $v \in \C^n$ we have the estimate
\begin{eqnarray*}
\frac{\parallel v \parallel}{2 L_G(p,v)} \le F_G(p,v) \le \frac{\parallel v \parallel}{ L_G(p,v)}
\end{eqnarray*}
\end{lemma}
This result implies many useful consequences. For example, $G$ becomes convex after a biholomorphic change of coordinates near a smooth strictly pseudoconvex boundary point. 
Lemma \ref{BedPin} then implies that $F_G(z,v) \ge C/ dist (z, b\Omega)$ for vectors $v$ which are transverse (say, orthogonal) to the holomorphic tangent space to $bG$ at $p$. This implies the classical fact that a smoothly bounded strictly pseudoconvex domain is  complete hyperbolic.

\section{Proof of the main theorem}   Our approach is based on \cite{CoSu} and employs the scaling method due to Frankel \cite{F1}. However, in difference with respect to \cite{CoSu}
we do not use general results of Frankel on convergence of dilated families. Our proof is self-contained and uses only Lemma \ref{BedPin}.

 Assume that we are in hypothesis of Theorem \ref{theo1}.

\subsection{Scaling} Suppose that $\Omega$ is of the form (\ref{BoundPoint}) in a coordinate neighborhood $U$ of $p$. Recall that  the strictly pseudoconvex hypersurfaces $\Gamma_j = \{\rho_j = 0\}$  are called  faces of $b\Omega \cap U$.

\begin{lemma}
\label{normalization}
There exists a local biholomorphic change of coordinates $\Phi$ such that $\Phi(p) = 0$ and $\Phi(\Omega \cap U)$
is convex.
\end{lemma}
For the proof see Proposition 1.1 in \cite{CoSu}. In what follows we assume that the local coordinates are fixed accrding Lemma \ref{normalization} so that $\Omega \cap \varepsilon \B^n$ is convex for $\varepsilon > 0$ small enough (with some abuse of notations we identify $\Phi(\Omega \cap U)$ with $\Omega \cap \varepsilon \B^n$). Note that in these coordinates
every local defining function of $\Omega$ near the origin has the expansion
\begin{eqnarray}
\label{exp1}
\rho_j(z) = \Re z_j + H_j(z,\overline z) + S_j(z ), j= 1,...,m
\end{eqnarray}
where each $H_j$ is a positive defined Hermitian quadratic form and $S _j(z) = o(\vert z \vert^2)$.

Let $\Omega_k = (f^k)^{-1}(\Omega \cap \varepsilon \B^n)$. Since the sequence $(f^k)$ converge to $0$ uniformly on compact subsets of $\Omega$, one can assume that 
 $(\Omega_k)_k$ is an increasing sequence of domains in $\Omega$ such that $\Omega = \cup_k \Omega_k$. 

Fix a point $q$ which belong to $\Omega_k$ for all $k$. Set $p^k := f^k(q)$ and consider the affine linear maps $$A^k(z) := (df^k(q))^{-1}(z - p^k).$$ Define a new sequence of maps 
\begin{eqnarray}
\label{S0}
g^k:= A^k \circ f^k.
\end{eqnarray} Note that 
\begin{eqnarray}
\label{S1}
g^k(q) = 0 \, \, \mbox{ and }  dg^k(q) = Id \,\, \mbox{for all} \,\, k
\end{eqnarray}
Consider the images $G_k = g^k(\Omega_k) = A^k(\Omega \cap \varepsilon \B^n)$. Our ultimate goal is to prove that the sequence of convex domains $(G^k)$ converges in the Hausdorff distance to a domain $G$ and to determine this limit domain $G$.
 
\subsection{Convergence of domains} First we note that the tangent maps $$R^k:= df^k(q)$$ converge to $0$; therefore, the domains $(df^k(q))^{-1}(\varepsilon \B^n)$ converge to the whole space $\C^n$. For this reason they do not affect our argument and we do not mention them anymore. Every domain $G_k$ is defined by 
\begin{eqnarray*}
\{ z : \rho_j(p^k + R^k z) < 0, j= 1,...,m \}
\end{eqnarray*}
Set $\tau_k^j:= \vert \rho_j(p^k) \vert$ and $\delta_k:= \inf_j \tau^j_k$.  Consider the functions
\begin{eqnarray*}
\phi_j^k(z) = (\tau_k^j)^{-1} \rho_j(p^k + R^k z), j= 1,...,m.
\end{eqnarray*}
Their expansions at the origin have the from
\begin{eqnarray}
\label{exp2}
\phi_j^k(z) = -1 + \Re \lambda_j^k(z) +  (\tau_k^j)^{-1}\Re Q_j^k(R^k z, R^k z) + (\tau_k^j)^{-1}H^k_j(R^k z, \overline{R^k z}) +  S_j^k(z)
\end{eqnarray}
Here $\lambda_j^k$ are complex linear forms, $Q_j^k(w,w)$ are holomorphic quadratic forms and in view of (\ref{exp1}) one has 
$Q_j^k \to 0$ as $k \to \infty$; $H^k_j(w, \overline w)$ are positive defined quadratic forms converging respectively to $H_j$ from (\ref{exp1}).
Finally, $S^k_j(z) = o(\vert z \vert^2)$ uniformly in $k$.
\begin{lemma}
\label{convergence1}
For every $j$, the sequence $(\phi_j^k)_k$ converges (after passing to a subsequence) uniformly on compact subsets of $\C^n$ as $k \to \infty$, to the function
$$\phi_j = -1 + Re \lambda_j(z) + H_j'(z,\overline z)$$
Here every $\lambda_j$ is a complex linear from and every $H_j'$ is a non-negative Hermitian quadratic from.
\end{lemma}
\proof There exists $C > 0$ such that for every $k$ and $j$ 
$$C^{-1} dist (p^k,\Gamma_j) \le \tau_j^k \le C dist(p^k, \Gamma_j)$$
Since $dist(p^k,b\Omega) = \inf_j dist(p^k,\Gamma_j)$, one has
$$C^{-1} dist (p^k, b\Omega) \le \delta_k \le C dist(p^k, b\Omega)$$
Lemma \ref{belowbound} and Lemma \ref{upperbound} imply the estimates (for each $v \in \C^n$):
\begin{eqnarray*}
C^{-1} \parallel v \parallel \ge F_{\Omega_k}(q,v) \ge F_{\Omega \cap \varepsilon \B^n} \ge \frac{C \parallel v \parallel}{\delta_k^{1/2}} \ge \frac{C \parallel v \parallel}{(\tau_k^j)^{1/2}}
\end{eqnarray*}
which gives
\begin{eqnarray}
\label{est1}
\parallel R^k \parallel \le C (\tau_k^j)^{1/2}
\end{eqnarray}
As a consequence we obtain that in (\ref{exp2}) the sequence $(\tau_k^j)^{-1}\Re Q_j^k(R^k z, R^k z)$ converges to $0$ uniformly on compact subset of $\C^n$ and 
$(\tau_k^j)^{-1}H^k_j(R^k z, \overline{R^k z}) $ converges uniformly on compact subsets of $\C^n$. It is also easy to see that $S^k_j$ converges to $0$ uniformly on compact subsets of $\C^n$.

 Next, it follows by (\ref{S1}) that 
$F_{\Omega_k}(q,v) = F_{G_k}(0,v)$ for all $k$ anf $v \in \C^n$. Therefore, there exists $C > 0$ such that 
$$C^{-1} \parallel v \parallel \le F_{G_k}(0,v) \le C \parallel v \parallel$$
Since the domains $G_k = \{ \phi_j^k < 0, j=1,...,m \}$ are convex,  by Lemma \ref{BedPin} one has
\begin{eqnarray}
\label{est3}
C^{-1} \le L_{G_k}(0,v) \le C
\end{eqnarray}
for all $k$ and $v$. 
Arguing by absurd, assume that the norms $\alpha_j^k$ of the forms $\lambda^k_j$ are not bounded in $k$; one can assume that $\alpha_j^k \to \infty$.  Then the functions $(\alpha_j^k)^{-1}\phi_j^k$ converge to functions $Re \theta_j(z)$, where $\theta_j$ is a non-zero complex linear form. This means that the boundaries of convex domains $G_k$
approach the origin as $k \to \infty$ and for some non-zero vector $v$ one has $L_{G_k}(0,v) \to 0$ as $k \to \infty$.This contradiction proves that the sequence of norms of the forms $\lambda_j^k$ is bounded and concludes the proof of Lemma.

Thus, the domains $G_k$ converge in the Hausdorff distance to the domain

\begin{eqnarray}
\label{limit}
G = \{ z: \phi_j(z) < 0, j = 1,...,m \}
\end{eqnarray}
Our goal now is to prove that $m= 1$ and $G$ is biholomorphic to $\B^n$.

\subsection{Identification of $G$. Case $m=1$:   the reheating}
First we consider the simplest case whee $m = 1$.  Then 
$$G = \{ z: -1 + Re \lambda(z) + H(z,\overline z) < 0 \}$$
If a non-zero vector $v$ is contained in the intersection  $\ker \lambda \cap \ker H = \{ 0 \}$, then the complex line through the origin in the direction of $v$ is contained in $G$ and  $L_{G}(0,v) = \infty$ which contradicts to (\ref{est3}).
Hence the restriction of $H$ on $\ker \lambda$ is positive defined and $G$ is biholomorphic to $\B^n$.

In order to conclude the proof of Theorem in this case  we need the following

\begin{lemma}
\label{MapConvergence}
The sequence of maps (\ref{S0}) converges (after passing to a subsequence) to a biholomorphism between $\Omega$ and $G$.
\end{lemma}
\proof Fix a compact subset $K \in G$. Since the sequence of convex domains $(G_k$) converges to $G$,  there exists $k_0$ such that $K \subset G_k$ for all $k \ge k_0$. It follows from Lemma \ref{BedPin} that there exists $C > 0$ such that 
\begin{eqnarray*}
F_{G_k}(z,v) \ge C\parallel v \parallel \,\, \mbox{for all} \,\, z \in K, \,\, v \in \C^n \,\, \mbox{and} \,\, k \ge k_0.
\end{eqnarray*} Then the classical argument (see \cite{R}) shows that the family $(g^k)$ is normal. Since $g^k(q) = 0$ for all $k$, the sequence  $(g^k)$ contains a subsequence converging uniformly on compact subsets of $\Omega$ to a holomorphic map $g$. On the other hand, the domain  $\Omega$ is hyperbolic (Corollary \ref{Hyp}).  Hence a similar  argument 
implies the convergence of the family of inverse maps $((g^k)^{-1})$. Now the classical theorem of H.Cartan (see \cite{N}) shows that $g:\Omega \to G$ is biholomorphic.

\subsection{Identification of $G$. Case $m>  1$:  the general case }
We consider now the case where $m > 1$. In this case it is appropriate to modify slightly the scaling sequence $(g^k)$. Namely, consider the linear mappings
$$B^k: z \mapsto w$$
$$w_j = \sum_{l=1}^n \frac{\partial \rho_j}{\partial z_l}(p_k) z_l, j=1,...,m,$$
$$w_j = z_j, j= m+1,...,n$$
Note that the sequence of linear maps $(B^k)_k$ converges to the identity map as $k$ ends to $\infty$. Consider the sequence of maps 
\begin{eqnarray*}
\tilde  g^k:= (R^k)^{-1} \circ B^k \circ (f^k - p^k)
\end{eqnarray*}
and consider the domains $$\tilde G_k = \tilde g^k(\Omega_k \cap U).$$ Then 
$$\tilde G^k = \{ z: \rho_j(p^k + (B^k)^{-1} \circ R^k z) < 0, j= 1,...m \}$$
Precisely  as above,  consider the functions

\begin{eqnarray*}
\phi_j^k(z) = (\tau_k^j)^{-1} \rho(p^k + (B^k)^{-1} \circ R^k z)
\end{eqnarray*}
Suppose that $$R^k z = (r^k_1 z,..., r^k_n z)$$ where $r^k_j$ are complex linear forms. Then the expansions at the origin have the from
\begin{eqnarray}
\label{exp3}
\phi_j^k(z) = -1 + (\tau_k^j)^{-1} \Re r_j^k(z) +  (\tau_k^j)^{-1}\Re Q_j^k(R^k z, R^k z) + (\tau_k^j)^{-1}H^k_j(R^k z, \overline{R^k z}) +  S_j^k(z)
\end{eqnarray}
Here $r_j^k$ are complex linear forms as above (the components of $R^k$; the additional maps $B^k$ are introduced to the scaling process in order to make 
appear $r^k_j$ explicitely in this expansion), $Q_j^k(w,w)$ are holomorphic quadratic forms and in view of (\ref{exp1}) one has 
$Q_j^k \to 0$ as $k \to \infty$; $H^k_j(w, \overline w)$ are positive defined quadratic forms converging respectively to $H_j$ from (\ref{exp1}).
Finally, $S^k_j(z) = o(\vert z \vert^2)$. Here we have used the fact that the sequence $(B^k)$ converges to identity.

Now Lemma \ref{convergence1} can be applied for every $j$. We obtain that $(\phi_j^k)_k$ converges (after passsing to a subsequence) 
to 
$$\phi_j = -1 + \Re \lambda_j(z) + H_j'(z,\overline z), j= 1,...m$$
Here every $\lambda_j$ is a complex linear form and every $H_j'$ is a  Hermitian quadratic from. Notice that the forms $H_j’$ are non-negative, hence the functions $\phi_j$ are plurisubharmonic.
The key observation is that we can obtain more information about the limit functions.
Namely, since the norms of linear forms $(\tau_k^j)^{-1} r_j^k(z)$ are bounded, the norms of the forms $(\tau_k^j)^{-1/2} r_j^k(z)$ tend to $0$. Therfore, the functions $\phi_j$
 can be written as $$\phi_j(z) = -1 + \Re\lambda_j(z) + \tilde H_j(L_{m+1}(z),...,L_n(z), \overline{L}_{m+1}(z),...,\overline{L}_n(z))$$
where $\lambda_j$, $L_j$, $j =1,...n$ are complex linear forms. Furthermore, these forms are linearly independent. Indeed, if a non-zero vector $v$ is contained in the intersection 
of their kernels, then the complex line through the origin in the direction of $v$ is contained in the limit domain $G = \{ \phi_j < 0, j= 1,...,m \}$. But this contradicts  the hyperbolicity of the convex domain $G$ at the origin, as above. Hence, after a complex linear change of coordinates $G$ becomes
\begin{eqnarray}
\label{image1}
G = \{ z: \psi_j = \Re s_j + F_j(t, \overline{t}) < 0, j= 1,...,m \}
\end{eqnarray}
where we put $s = (z_1,...,z_m)$ and $t = (z_{m+1},...,z_{n})$,  and every Hermitian from $F_j$ is positive defined on the space $\C^{n-m}(t)$.
Then, as it is easy to see, the limit domain $G$ is hyperbolic . In fact, $G$ is biholomorphic to a bounded domain 
\begin{eqnarray}
\label{image2}
\tilde G = \{ z: \tilde\psi_j = \vert s_j \vert^2 + \tilde F_j(t, \overline{t}) < 0, j= 1,...,m \}
\end{eqnarray}
where the forms $\tilde F_j$ are positive defined on $\C^{m}(t)$. In order to see this fact,  it suffices to apply the Caley transform to every defining function
$\psi_j$ in the space $\C^{m+1}(s_j,t)$.

As in Lemma \ref{MapConvergence}, now we conclude that  the family $(\tilde g^k)$ is normal. Hence, $\Omega$ is biholomorphic to $G$ as above.

The last Step of the proof is to show that $m=1$.

Let $f: \tilde G \to \Omega$ be a biholomorphic map. Let  $a$ be a boundary point of $\tilde G$ contained in the claster set of $f^{-1}$ at $p$. This means that  there exists a sequence $(a^k)$
in $\tilde G$ converging to $a$ such that $f(a^k)$ converges to $p$.
Then the map $f$ extends as a Holder continuous map on the boundary $b{\tilde G}$ in a neighborhood of $a$ by Theorem 1.1 of \cite{Su}. Note that in \cite{Su} a boundary point from the source domain ($\tilde G$ in our case) is required to be piecewise smooth strictly pseudoconvex. However, this assumption is imposed there because  maps under consideration in \cite{Su} 
are only locally proper. In our case $f$ is biholomorphic, the domain $\tilde G$  admits a global defining plurisubharmonic function $\sup_j \tilde \phi_j$ and Step 2 of the proof \cite{Su} (based on the Hopf lemma) goes through directly. The remaining part of the proof from \cite{Su} goes through literally which establishes the Holder continuity up to the boundary. 

Now, when $m > 1$, every face of $b\tilde G$ is foliated by complex lines. Then the argument of  Theorem 1.2 of \cite{Su} or \cite{CoSu} shows that the Jacobian determinant of $f$ vanishes identically on a one-sided neighborhood of $a$ in $\tilde G$ and, therefore, everywhere on $\tilde G$. This is a contradiction because $f$ is a biholomorphic map. Hence, $m = 1$. The proof is finished.

\bigskip

\bigskip

{\small

 Universit\'e  de Lille, Laboratoire
Paul Painlev\'e,
U.F.R. de
Math\'ematiques, 59655 Villeneuve d'Ascq, Cedex, France, e-mail: sukhov@math.univ-lille1.fr,
 and  Institut of Mathematics with Computing Centre , Ufa Federal Research Centre of Russian
Academy of Sciences, 45008, Chernyshevsky Str. 112, Ufa, Russia.

The author is partially suported by Labex CEMPI.
}

\end{document}